\documentclass[11pt,leqno]{article}
\usepackage{latexsym}
\usepackage{amsmath,amssymb}
\usepackage{amsthm}
\usepackage[spanish,english]{babel}
\title{Injectors in $\pi$-separable groups}
\author{M. Arroyo-Jord\'a, P. Arroyo-Jord\'a, R. Dark,\\ A.~D. Feldman, and M.~D. P\'{e}rez-Ramos}

\date{}
\newtheorem{teo}{Theorem}[section]
\newtheorem{lem}[teo]{Lemma}
\newtheorem{pro}[teo]{Proposition}
\newtheorem{cor}[teo]{Corollary}
\theoremstyle{definition}
\newtheorem*{remark}{Remark}
\newtheorem*{remarks}{Remarks}

\newtheorem{de}[teo]{Definition}

\DeclareMathOperator{\Syl}{Syl}
\DeclareMathOperator{\Hall}{Hall}

\DeclareMathOperator{\Core}{Core}
\DeclareMathOperator{\proj}{Proj}
\DeclareMathOperator{\cov}{Cov}
\DeclareMathOperator{\inj}{Inj}
\DeclareMathOperator{\G}{\overline{\it G}}
\DeclareMathOperator{\Q}{\overline{\it Q}}
\DeclareMathOperator{\P^*}{\overline{\it P^*}}
\DeclareMathOperator{\N}{\overline{\it N_G(V)}}
\DeclareMathOperator{\M}{\overline{\it M}}
\DeclareMathOperator{\V}{\overline{\it V}}
\DeclareMathOperator{\Vpi}{\overline{\it V}_{\pi}}
\DeclareMathOperator{\Vpin}{\overline{\it V}_{\pi^\prime}}
\DeclareMathOperator{\Ls}{\overline{\it L}}
\DeclareMathOperator{\R^*}{\overline{\it R^*}}
\DeclareMathOperator{\npi}{\mathfrak N^\pi}
\DeclareMathOperator{\F}{\mathfrak F}
\DeclareMathOperator{\p}{\mathbb{P}}

\begin{document}
\maketitle

\begin{abstract}

Let $\pi$ be a set of primes. We show that $\pi$-separable groups have a conjugacy class of $\mathfrak F$-injectors for suitable Fitting classes $\mathfrak F$, which coincide with the usual ones  when specializing to soluble groups.
\medskip

\medskip

\end{abstract}

\section{Introduction and preliminaries}

All groups considered are finite.

One of the fundamentals facts in the theory of finite soluble groups is the theorem of B. Fischer, W. Gasch\"{u}tz and B. Hartley,  which states the existence and conjugacy of $\F$-injectors in finite soluble groups for Fitting classes $\F$ (\cite{FGH}, \cite[VIII. Theorem (2.9), IX. Theorem (1.4)]{DH}). An important stream of research has considered the extent of the validity of this result to all finite groups. We refer to \cite{DH,BE} for background on classes of groups and for accounts about the development of this topic;  we shall adhere to their notations. We recall that if $\mathfrak F$ is a class of groups, an $\mathfrak F$-injector of a group $G$ is a subgroup $V$ of $G$ with the property that $V\cap K$ is an $\mathfrak F$-maximal subgroup of $K$ (i.e. maximal as subgroup of $K$ in $\mathfrak F$) for all subnormal subgroups $K$ of $G$. The existence of $\F$-injectors in all groups implies that $\F$ is a Fitting class, which is defined as a non-empty class of groups closed under taking normal subgroups and products of normal subgroups. The existence and properties of Carter subgroups, i.e. self-normalizing nilpotent  subgroups, in soluble groups are the cornerstone in the proof of the above-mentioned theorem of Fischer, Gasch\"{u}tz and  Hartley. In a previous paper~\cite{AADFP} we initiated the study of an extension of the theory of soluble groups to the universe of $\pi$-separable groups, $\pi$ a set of primes. We analyzed the reach of $\pi$-separability further from soluble groups, by means  of complement and Sylow bases and Hall systems, based on the remarkable property that $\pi$-separable groups have Hall $\pi$-subgroups, and every $\pi$-subgroup is contained in a conjugate of any Hall $\pi$-subgroup. We also proved that  $\pi$-separable groups have a conjugacy class of subgroups which specialize to Carter subgroups within the universe of soluble groups. The main results in the present paper, namely Lemma~\ref{lema*}, Theorem~\ref{teo*} and Corollary~\ref{cor*} below, show that these  Carter-like subgroups enable an extension of the existence and conjugacy of injectors to $\pi$-separable groups.

If $\pi$ is a set  of primes, let us recall that a group $G$ is $\pi$-separable if every composition factor of $G$ is either a $\pi$-group or a $\pi'$-group, where $\pi'$ stands for the complement of $\pi$ in the set $\p$ of all prime numbers. Clearly, $\pi$-separability is equivalent to
$\pi'$-separability, so that by the Feit-Thompson theorem, every $\pi$-separable group is either $\pi$-soluble or $\pi'$-soluble, where for any set of primes $\rho$, a group is $\rho$-soluble if it is $\rho$-separable with every $\rho$-composition factor a $p$-group for some prime $p\in \rho$. Also  Burnside's $p^aq^b$-theorem implies that $\pi$-separable groups are $\pi$-soluble if $|\pi|\le 2$. We refer to \cite{Go} for basic results on $\pi$-separable groups.

These remarks about $\rho$-solubility clarify the reach of our main results, which apply to $\pi'$-soluble groups once the set of primes $\pi$ is fixed.

More precisely, our study of $\pi$-separable groups relies on convenient extensions of the class $\mathfrak N$ of nilpotent groups as well as of normal subgroups, according with the set of primes $\pi$, as follows.

Let $\pi$ be a set of primes. Let $$\mathfrak N^\pi=\mathfrak  E_\pi \times \mathfrak  N_{\pi'}=(G=H\times K\mid H\in \mathfrak  E_\pi,\ K\in \mathfrak  N_{\pi'}),$$ where $\mathfrak  E_\pi$  denotes the class of all $\pi$-groups and $\mathfrak  N_{\pi'}$ the class of all nilpotent $\pi'$-groups.
In the particular cases when either $\pi=\emptyset$ or $\pi=\{p\}$, $p$ a prime, ($|\pi|\le 1$), then $\mathfrak N^\pi=\mathfrak  N$ is the class of all nilpotent groups.

We observe that $\mathfrak N^\pi$ is a saturated formation and appeal to the concept of $\mathfrak N^\pi$-Dnormal subgroup. We refer to \cite{AP, AADFP} for the concept of $\mathfrak G$-Dnormal subgroup for general saturated formations $\mathfrak G$, and specialize this definition to our particular saturated formation $\mathfrak N^\pi$ next. For notation, if $\rho$ is a set of primes and  $G$ is a group, $\Hall_{\rho}(G)$ denotes the set of all Hall $\rho$-subgroups of $G$. If $p$ is a prime, then $\Syl_p(G)$ stands for the set of all Sylow $p$-subgroups of $G$. If $G_\rho \in \Hall_\rho(G)$ and $H\le G$, we write $G_\rho\searrow H$ to mean that $G_\rho$ reduces into $H$, i.e. $G_\rho\cap H\in \Hall_\rho (H)$.

\begin{de}\textup{(\cite[Definition 3.1]{AP}, \cite{AADFP})} A subgroup $H$ of a group $G$ is said to be \emph{$\npi$-Dnormal} in $G$ if it satisfies the following conditions:
\begin{description}
\item[$(1)$] whenever $p\in \pi'$ and $G_p\in \Syl_p(G)$, $G_p\searrow H$, then $G_p\le N_G(H)$;
\item[$(2)$] whenever $p\in \pi$ and $G_p\in \Syl_p(G)$, $G_p\searrow H$, then
\begin{itemize}
\item $G_p\le N_G(H)\ $ if $\pi=\{p\}$,\ \ or
\item $G_p\le N_G(O^\pi(H))\ $ if $|\pi|\ge 2$.
\end{itemize}
\end{description}
\end{de}

Hence, for $\mathfrak N^\pi=\mathfrak  N$, $\mathfrak  N$-Dnormal subgroups are exactly normal subgroups. Note also that normal subgroups are $\mathfrak N^\pi$-Dnormal for any set of primes $\pi$.
\smallskip

 $\npi$-Dnormal subgroups are nicely characterized as follows.

\begin{pro}\label{pro1}\textup{\cite[Proposition 2.3]{AADFP}} Let $H$ be a subgroup  of a group $G$. Then:
\begin{enumerate}\item Assume that $|\pi|\le 1$. Then $\mathfrak N^\pi =\mathfrak  N$ and $H$ is $\mathfrak  N$-Dnormal in $G$ if and only if $H$ is normal in $G$.
\item Assume that $|\pi|\ge 2$. Then the following statements are equivalent:
\begin{itemize} \item[(i)] $H$ is $\mathfrak N^\pi$-Dnormal in $G$;
\item[(ii)] $O^{\pi}(H)\unlhd G$ and $O^{\pi}(G)\le N_G(H)$.
\end{itemize}
\end{enumerate}
\end{pro}

For our study of Fitting classes it is useful to introduce a corresponding extension  of subnormality.

\begin{de} A subgroup $S$ of a group $G$ is said to be \emph{$\npi$-Dsubnormal} in $G$ if there is a chain of subgroups $$S=S_0\le S_1\le \cdots\le S_k=G,$$
such that $S_i$ is $\npi$-Dnormal in $S_{i+1}$ if $0\le i\le k-1$.
\end{de}

As for $\npi$-Dnormal subgroups, if $\mathfrak N^\pi=\mathfrak  N$, then $\mathfrak  N$-Dsubnormal subgroups are exactly subnormal subgroups;  also  subnormal subgroups are $\mathfrak N^\pi$-Dsubnormal for any set of primes $\pi$.
\smallskip

In order to prove our main results we will need some properties of $\npi$-Dnormal and $\npi$-Dsubnormal subgroups that we gather in the next lemmas.

\begin{lem}\label{lem2} Let $H$ be a subgroup of a group $G$, $N\unlhd G$ and $g\in G$. Then:
\begin{enumerate}
\item If $H$ is $\npi$-Dnormal in $G$, then $H^g$ is $\npi$-Dnormal in $G$.
\item If $H$ is $\npi$-Dnormal in $G$ and $H\le L\le G$, then $H$ is $\npi$-Dnormal in $L$.
\item If $H$ is $\npi$-Dnormal in $G$, then $HN/N$ is $\npi$-Dnormal in $G/N$.
\item If $N\le H$ and $H/N$ is $\mathfrak N^\pi$-Dnormal in $G/N$, then $H$ is $\mathfrak N^\pi$-Dnormal in $G$.
\item If $N\le H$ and $H/N$ is $\mathfrak N^\pi$-Dsubnormal in $G/N$, then $H$ is $\mathfrak N^\pi$-Dsubnormal in $G$.
\item If $G\in \mathfrak N^\pi$, then $H$ is $\mathfrak N^\pi$-Dsubnormal in $G$.
\end{enumerate}
\end{lem}
\smallskip

{\it Proof.} If $|\pi|\le 1$, the result refers to normal and subnormal subgroups and it is clear. Assume that $|\pi|\ge 2$. Then:
\smallskip

1 and 2 are easily proven.
\smallskip

3. By Proposition~\ref{pro1}, assuming that $O^\pi(H)\unlhd G$ and $O^\pi(G)\le N_G(H)$, we need to prove that $O^\pi(HN/N)\unlhd G/N$ and $O^\pi(G/N)\le N_G(HN/N)$. But this is clear since for any $L\le G$, $O^\pi(LN/N)=O^\pi(L)N/N$.
\smallskip

4. Again by Proposition~\ref{pro1}, we assume that $O^\pi(H/N) = O^\pi(H)N/N  \unlhd G/N$ and $O^\pi(G/N) = O^\pi(G)N/N\le N_G(H/N),$ and need to prove that $O^\pi(H)\unlhd G$ and $O^\pi(G)\le N_G(H)$. The second property follows clearly, and also that $[G,O^\pi(H)]\le O^\pi(H)N\le H$. Then $[G, O^\pi(H)]\le H$, which is equivalent to $[G, O^\pi(H)]\le O^\pi(H)$ (see \cite[Remark 2.2]{AADFP}). Therefore, $G\le N_{G}(O^\pi(H))$, and we are done.
\smallskip

5. It follows from the definition of $\npi$-Dsubnormal subgroup together with part 4.
\smallskip

6. If $G\in \mathfrak N^\pi$ then $G=G_\pi\times G_{\pi'}$, $G_\pi=O_\pi(G)$, and $G_{\pi'}=O_{\pi'}(G)\in \mathfrak N$. Also $H\in \npi$ so that
$H=H_\pi\times H_{\pi'}$, $H_\pi=O_\pi(H)$, and $H_{\pi'}=O_{\pi'}(H)\in \mathfrak N$. Then we deduce easily that $H$ is $\npi$-Dnormal in $G_\pi\times H_{\pi'}$. Since $G_{\pi'}$ is nilpotent, $H_{\pi'}$ is subnormal in $G_{\pi'}$. Hence $G_\pi\times H_{\pi'}$ is subnormal in $G$, and so also $\npi$-Dsubnormal. Therefore, $H$ is $\npi$-Dsubnormal in $G$, as desired.\qed

\begin{lem}\label{lem3} Let $H$ be an $\npi$-Dnormal subgroup of a group  $G$. Then:
\begin{enumerate}
\item $H/O^\pi(H)\le O_\pi(G/O^\pi(H)).$
\item Let $C=\Core_G(H)$ and $\langle H^G\rangle $ be the normal closure of $H$ in $G$. Then $\langle H^G\rangle /C\le O_\pi(G/C) $; equivalently, $\langle H^G\rangle /C\in \mathfrak  E_\pi$.
\item If $V\le G$, then $H\cap V$ is $\npi$-Dnormal in $V$.
\end{enumerate}
\end{lem}
\smallskip

{\it Proof.} 1. By Proposition~\ref{pro1} we know that $O^\pi(H)\unlhd G$, and we aim to prove that $H/O^\pi(H)\le O_\pi(G/O^\pi(H)).$ We argue by induction on $|G|$. If $O^\pi(H)\neq 1$, the result is clear by Lemma~\ref{lem2}(3) and inductive hypothesis. We may then assume that $O^\pi(H)= 1$, i.e. $H$ is a $\pi$-group, and we need to prove that $H\le O_\pi(G)$. If $O^\pi(G)=1$, then $H\le G=O_\pi(G)$, and we are done. Consider now the case when $O^\pi(G)\neq 1$. Note that $O^\pi(G)\le N_G(H)$ by Proposition~\ref{pro1}, because $H$ is $\npi$-Dnormal in $G$. Then  $H\cap O^\pi(G)\le O_\pi(O^\pi(G))\le O_\pi(G)$. If  $H\cap O^\pi(G)\neq 1$, again by Lemma~\ref{lem2}(3) and inductive hypotesis it follows that $HO_\pi(G)/O_\pi(G)\le O_\pi(G/O_\pi(G))=O_\pi(G)/O_\pi(G)$, and so $H\le O_\pi(G)$ as desired. If   $H\cap O^\pi(G)= 1$, then $[H,O^\pi(G)]=H\cap O^\pi(G)=1$ and so $H\le C_G(O^\pi(G))$. In the case when $C_G(O^\pi(G))<G$, by Lemma~\ref{lem2}(2) and inductive hypothesis we deduce that $H\le O_\pi(C_G(O^\pi(G)))\le O_\pi(G)$. Assume finally that $C_G(O^\pi(G))=G$. Then
$O^\pi(G)=O_{\pi'}(G)$ and, by the Schur-Zassenhaus theorem and the fact that $O^\pi(G)\le Z(G)$, there is a unique Hall $\pi$-subgroup of $G$, which is $O_\pi(G)$, and also $H\le O_\pi(G)$, which concludes the proof.
\smallskip

2. This is clear from part 1.
\smallskip

3. If $|\pi|\le 1$, then $H\unlhd G$ and the result is clear. Assume that $|\pi|\ge 2$. Then, by Proposition~\ref{pro1}, assuming that $O^\pi(G)\le N_G(H)$ and $O^\pi(H)\unlhd G$, we need to prove that $O^\pi(V)\le N_V(V\cap H)$ and $O^\pi(V\cap H)\unlhd V$. Since $O^\pi(V)\le O^\pi(G)$, it is clear that $O^\pi(V)\le N_V(V\cap H)$.  We consider $O^\pi(V\cap H)=\langle X\mid X\ \pi'\text{-subgroup of }V\cap H\rangle$. Let $x\in V$. If $X$ is a $\pi'$-subgroup of $V\cap H$, then $X^x\le V\cap O^\pi(H)^x=V\cap O^\pi(H)\le V\cap H$, which implies that $X^x\le O^\pi(V\cap H)$. Hence, $O^\pi(V\cap H)^x\le O^\pi(V\cap H)$, and so $O^\pi(V\cap H)\unlhd V$, which concludes the proof.\qed

\begin{lem}\label{maxFDn} If $M$ is a maximal $\npi$-Dnormal proper subgroup of a $\pi'$-soluble group $G$, then $G^{\npi}\le M$, where $G^{\npi}$ denotes the $\npi$-residual of $G$, i.e. the smallest normal subgroup in $G$ with quotient group an $\npi$-group.\end{lem}
\smallskip

{\it Proof.} By Lemma~\ref{lem3}(2) we know that $\langle M^G\rangle /\Core_G(M)\in \mathfrak E_\pi$. Since $M$ is a maximal $\npi$-Dnormal proper subgroup of $G$ we deduce that either $\langle M^G\rangle=M$, i.e. $M\unlhd G$, or $\langle M^G\rangle=G$. In the first case, since $G$ is $\pi'$-soluble, either $O^\pi(G)\le M$ or $O^p(G)\le M$ for some $p\in \pi'$, which imply that $G^{\npi}\le M$. If $\langle M^G\rangle=G$, then $G/\Core_G(M)\in \mathfrak E_\pi$, and so also $G^{\npi}\le M$.\qed

\section{Injectors in $\pi$-separable groups}

This section is devoted to proving our main results. We fix a set of primes $\pi$ and introduce suitable Fitting classes, which we will call $\npi$-Fitting classes, as defined below in this section. They appear to be Fitting classes  with stronger closure properties involving $\npi$-Dnormal subgroups. Then we prove the existence and conjugacy of associated injectors in  $\pi'$-soluble groups. (See Theorem~\ref{teo*}, Corollary~\ref{cor*}, below).

Our treatment adheres to the approach in \cite[Chapter VIII]{DH}, and we present our study within the framework of Fitting sets, instead of  general Fitting classes. As mentioned in \cite[VIII.1]{DH}, it does not cause any additional difficulty and has important advantages, both in terms of scope of results and working techniques.

We extend the concept of Fitting set \cite[VIII. Definition (2.1)]{DH} to $\mathfrak N^\pi$-Fitting set as follows, by replacing the terms ``normal subgroup'' and ``subnormal subgroup'' by ``$\mathfrak N^\pi$-Dnormal subgroup'' and  ``$\mathfrak N^\pi$-Dsubnormal subgroup'', respectively, in the original definition.

\begin{de}\label{defs} A non-empty set $\mathcal F$ of subgroups of a group $G$ is called an \emph{$\mathfrak N^\pi$-Fitting set} of $G$ if the following conditions are satisfied:
\begin{description}\item {\bf FS1:} If $T$ is an $\mathfrak N^\pi$-Dsubnormal subgroup of $S \in \mathcal F$, then $T\in \mathcal F$;

\item {\bf FS2:} If $S,T\in \mathcal F$ and $S,T$ are $\mathfrak N^\pi$-Dnormal subgroups in $\langle S,T\rangle$, then $\langle S,T\rangle\in \mathcal F$;

\item {\bf FS3:} If $S \in \mathcal F$ and $x\in G$, then $S^x\in \mathcal F$.
\end{description}
\end{de}

Since normal subgroups are $\mathfrak N^\pi$-Dnormal, so subnormal subgroups are $\mathfrak N^\pi$-Dsubnormal, it is clear that an $\mathfrak N^\pi$-Fitting set is a Fitting set. Also, Fitting sets are exactly $\mathfrak N$-Fitting sets, for $\npi=\mathfrak N$ with $|\pi|\le 1$.

For the basic results about Fitting sets we refer to \cite[VIII.2]{DH}.

We recall in particular that for a Fitting set $\cal F$ of a group $G$, the  $\mathcal F$-radical of $G$, denoted $G_{\cal F}$ is the join of all normal
 $\cal F$-subgroups of $G$; or equivalently, the join of all subnormal $\cal F$-subgroups of $G$. For a subgroup $H$ of $G$,  we set $H_{\cal F}$ for the radical of $H$ associated to its Fitting set $\mathcal F_H=\{S\le H\mid S\in \mathcal F\}$ which we shall denote simply as $\mathcal F$. (See \cite[VIII. Definitions (2.3), Proposition (2.4)]{DH}.)

The following property is often useful:

\begin{lem}\textup{\cite[VIII. Proposition (2.4)(b)]{DH}}\label{radical} Let $\cal F$ be a Fitting set of a group $G$, and let $H\le G$ and $x\in G$. Then $(H_{\mathcal F})^x=(H^x)_{\cal F}$. In particular, $N_G(H)\le N_G(H_{\mathcal F})$.
\end{lem}

Our first aim is to prove that if $\cal F$ is an $\npi$-Fitting set, the $\cal F$-radical is equally described as the join of all $\npi$-Dnormal
 $\cal F$-subgroups, and also as the join of all $\npi$-Dsubnormal $\cal F$-subgroups. (Proposition~\ref{cor1} below.) The next lemma supplies a basic fact about the join  of $\npi$-Dnormal subgroups.

\begin{lem}\label{Dn} If $H,K$ are $\npi$-Dnormal subgroups of a group $G$, then $\langle H,K\rangle$ is $\npi$-Dnormal in $G$.
\end{lem}
\smallskip

{\it Proof.} If $|\pi|\le 1$, then $\npi$-Dnormal subgroups are exactly normal subgroups and the result is clear. Assume that $|\pi|\ge 2$. We argue by induction on $|G|$. By Proposition~\ref{pro1} we have that $O^\pi(H)$ and $O^\pi(K)$ are normal subgroups of $G$. If either $O^\pi(H)\neq 1$ or $O^\pi(K)\neq 1$, then the result is easily deduced by Lemma~\ref{lem2}, parts (3), (4), and inductive hypothesis. Assume that $O^\pi(H)=O^\pi(K)= 1$. By Lemma~\ref{lem3}(1), $H,K\le O_\pi(G)$. Consequently, $\langle H,K\rangle$ is a $\pi$-group, and $O^\pi(\langle H,K\rangle)=1$ is a normal subgroup of $G$. On the other hand, by Proposition~\ref{pro1}, $O^\pi(G)\le N_G(H)\cap N_G(K)$, and so $O^\pi(G)\le N_G(\langle H,K\rangle)$. Again  Proposition~\ref{pro1} implies finally that $\langle H,K\rangle$ is $\npi$-Dnormal in $G$.\qed

\begin{lem}\label{Dsn} Let  $\mathcal F$ be an $\mathfrak N^{\pi}$-Fitting set of a group $G$. If $H$ is an $\npi$-Dnormal subgroup of $G$, then $H_{\cal F}$ is $\npi$-Dnormal in $G$.
\end{lem}
\smallskip

{\it Proof.} If $|\pi|\le 1$, then again $\npi$-Dnormal subgroups are exactly normal subgroups and the result follows from Lemma~\ref{radical}. Assume that $|\pi|\ge 2$. By Proposition~\ref{pro1} we need to prove that $O^\pi(G)\le N_G(H_{\cal F})$ and $O^\pi (H_{\cal F})\unlhd G$.

Since $H$ is $\npi$-Dnormal in $G$, again Proposition~\ref{pro1} and Lemma~\ref{radical} imply that $O^\pi(G)\le N_G(H_{\cal F})$. We claim now that $O^\pi(H_{\cal F})= O^\pi(O^\pi(H)_{\cal F})$. This will imply that $O^\pi(H_{\cal F})\unlhd G$ by Proposition~\ref{pro1} and Lemma~\ref{radical}, as above, and we will be done.

It is not difficult to check that $O^\pi(H_{\cal F})\le O^\pi(H)_{\cal F}\le H_{\cal F}$. Then $O^\pi(O^\pi(H)_{\cal F})\le O^\pi(H_{\cal F})$, $O^\pi(O^\pi(H)_{\cal F})\unlhd H_{\cal F}$ and $H_{\cal F}/O^\pi(O^\pi(H)_{\cal F})$ is a $\pi$-group. Hence  $O^\pi(H_{\cal F})= O^\pi(O^\pi(H)_{\cal F})$ and the claim is proven.\qed

\begin{pro}\label{cor1} If $\mathcal F$ is an $\mathfrak N^\pi$-Fitting set of a group $G$, then \begin{align*}G_{\mathcal F}&=\langle H\le G\mid H\unlhd G,\ H\in \mathcal F\rangle = \langle H\le G\mid H\vartriangleleft\vartriangleleft G,\ H\in \mathcal F\rangle\\
&=\langle H\le G\mid H\ \mathfrak N^\pi\mbox{-Dnormal in }G,\ H\in \mathcal F\rangle\\& = \langle H\le G\mid H\ \mathfrak N^\pi\mbox{-Dsubnormal in }G,\ H\in \mathcal F\rangle.\end{align*}
\end{pro}
\smallskip

{\it Proof.} Set $R= \langle H\le G\mid H\ \mathfrak N^\pi\mbox{-Dnormal in }G,\ H\in \mathcal F\rangle$. Since normal subgroups are $\npi$-Dnormal, it  is clear that $G_{\mathcal F}\le R$. Moreover, it follows from the definition of $\mathfrak N^\pi$-Fitting set and Lemmas~\ref{lem2}(1) and \ref{Dn}, that $R$ is a normal subgroup of $G$ in $\cal F$, which implies that $R\le G_{\mathcal F}$. Consequently, $G_{\mathcal F}= R$.

Let now $S=\langle H\le G\mid H\ \mathfrak N^\pi\mbox{-Dsubnormal in }G,\ H\in \mathcal F\rangle$. It is clear that $G_{\mathcal F}\le S$. Let $H$ be an $\npi$-Dsubnormal subgroup of $G$. We claim that $H_{\cal F}\le G_{\mathcal F}$. Hence, if in addition $H\in \cal F$, then $H\le G_{\mathcal F}$. It will follow that $S\le G_{\mathcal F}$, which will conclude the proof.

If $H$ is an $\npi$-Dsubnormal subgroup of $G$, there is a chain of subgroups $H=H_0\le H_1\le \cdots\le H_k=G,$
such that $H_i$ is $\npi$-Dnormal in $H_{i+1}$ if $0\le i\le k-1$. For each $i=0,\ldots, k-1$, Lemma~\ref{Dsn} implies that $(H_i)_{\cal F}$ is $\npi$-Dnormal in $H_{i+1}$, and then $(H_i)_{\cal F}\le (H_{i+1})_{\cal F}$. Therefore, $H_{\cal F}=(H_0)_{\cal F}\le (H_1)_{\cal F}\le \cdots\le (H_k)_{\cal F}=G_{\mathcal F},$ which proves the claim.\qed

\begin{remark} As it might be expected, as a consequence of Proposition~\ref{cor1}, in the definition of $\npi$-Fitting set (Definition~\ref{defs}), $\npi$-Dnormal subgroups can be equivalently replaced by  $\npi$-Dsubnormal subgroups in condition (FS2).
\end{remark}

\begin{itemize}
\item Let $\cal F$ be an $\npi$-Fitting set of a  group $G$. Then:
\smallskip

If $S,T\in \mathcal F$ and $S,T$ are $\mathfrak N^\pi$-Dsubnormal subgroups in $\langle S,T\rangle$, then $\langle S,T\rangle\in \mathcal F$

\smallskip

{\it Proof.} By Proposition~\ref{cor1} we deduce that $S,T\le \langle S,T\rangle _{\cal F}$, and then $\langle S,T\rangle =\langle S,T\rangle _{\cal F}\in \cal F$.\qed
\end{itemize}

We recall also the concept of $\cal X$-injector of a group for a set of subgroups $\cal X$ of the group.

\begin{de}\textup{\cite[VIII. Definition (2.5)(b)]{DH}} Let $\cal X$ be a set of subgroups of a group $G$. An \emph{$\cal X$-injector} of $G$ is a subgroup $V$ of $G$ with the property that $V\cap K$ is an $\cal X$-maximal subgroup of $K$ (i.e. maximal as subgroup of $K$ in $\cal X$) for every subnormal subgroup $K$ of $G$. We shall denote the (possibly empty) set of $\cal X$-injectors of $G$ by $\inj_{\mathcal X}(G)$.
\end{de}

The following results about existence and properties of $\mathfrak N^\pi$-projectors in $\pi'$-soluble groups are essential for our purposes.

Let us recall that given a class of groups $\mathfrak  X$, a subgroup $U$ of a group $G$ is called an $\mathfrak  X$-projector of $G$ if $UK/K$ is an $\mathfrak  X$-maximal subgroup of $G/K$   for all $K\unlhd G$. The (possibly empty) set of $\mathfrak  X$-projectors of $G$ will be denoted by $\proj_{\mathfrak  X}(G)$.

Also, an $\mathfrak  X$-covering subgroup of $G$ is a subgroup $E$ of $G$ with the property that $E\in \proj_{\mathfrak  X}(H)$ whenever $E\le H\le G$. The set of $\mathfrak  X$-covering subgroups of $G$ will be denoted by $\cov_{\mathfrak  X}(G)$.

We refer to \cite{DH,BE} for convenient background about projectors and covering subgroups.

As in \cite[Definition 3.6]{AADFP}, we say that a subgroup $H$ of a group $G$ is \emph{self-$\mathfrak N^\pi$-Dnormalizing} in $G$ if whenever $H\le K\le G$ and $H$ is $\mathfrak N^\pi$-Dnormal in $K$, then $H=K$.

\begin{lem}\textup{(\cite[Theorem 3.4]{AADFP})}\label{teo1} For all $\pi'$-soluble groups $G$, $\emptyset \neq \proj_{\mathfrak N^\pi}(G)= \cov_{\mathfrak N^\pi}(G)$ and it is a conjugacy class of $G$.
\end{lem}

\begin{lem}\textup{(\cite[Theorem 3.13]{AADFP})}\label{teo2} For a subgroup $H$ of a $\pi'$-soluble group $G$ the following statements are pairwise equivalent:
\begin{enumerate}\item $H$ is an $\mathfrak N^\pi$-projector of $G$.
\item $H$ is an $\mathfrak N^\pi$-covering subgroup of $G$.
\item $H\in \mathfrak N^\pi$ is a self-$\mathfrak N^\pi$-Dnormalizing subgroup of $G$ and $H$ satisfies the following property:
\begin{gather}\text{If}\  H\le X\le G,\ \text{then}\ H\cap X^{\mathfrak N^\pi}\le (X^{\mathfrak N^\pi})'.\tag{*}
\end{gather}

\end{enumerate}
\end{lem}

The next result extends Hartley's result \cite[VIII. Lemma (2.8)]{DH} for $\mathfrak N^\pi$-Fitting sets and $\pi'$-soluble groups.

\begin{lem}\label{lema*} Let $\mathcal F$ be an $\mathfrak N^\pi$-Fitting set of a $\pi'$-soluble group $G$. Let $K$ be a normal subgroup of $G$ containing the $\mathfrak N^\pi$-residual $G^{\mathfrak N^\pi}$ of $G$, let $W$ be an $\mathcal F$-maximal subgroup of $K$, and let $V$ and $V_1$ be $\mathcal F$-maximal subgroups of $G$ which contain $W$.
\begin{description}\item (a) If $W\unlhd K$, then $V=(WP)_{\mathcal F}$, where $P$ is a suitable $\mathfrak N^\pi$-projector of $G$.
\item (b) In any case $V$ and $V_1$ are conjugate in $\langle V,V_1\rangle$. More precisely, there exists $x\in \langle V,V_1\rangle^{\mathfrak N^\pi}$ such that $V_1^x= V$.
\end{description}
\end{lem}
\smallskip

{\it Proof.} We mimic the proof of \cite[VIII. Lemma (2.8)]{DH}. The arguments get a bit more involved mainly by the fact that  $\mathfrak N^\pi$-projectors are not characterized as  self-$\mathfrak N^\pi$-Dnormalizing $\mathfrak N^\pi$-subgroups, in order to play the role of Carter subgroups in the original proof, but the characterization of $\mathfrak N^\pi$-projectors in Lemma~\ref{teo2} can be instead used to prove the result.
\smallskip

(a) We argue by induction on $|G|$. If $W\unlhd K$, then $W=K_{\mathcal F}$, and so also $W\unlhd G$ by Lemma~\ref{radical}.
We gather the  following facts which will be useful in the proof, where $L$ denotes any subgroup of $G$ containing $W$ and $U$ denotes  an $\mathcal F$-maximal subgroup of $L$ containing $W$:
\begin{enumerate}

 \item $L$ satisfies the hypotheses of the statement by considering $K\cap L$ and $U$ playing the role of $K$ and $V$, respectively.
    \smallskip

    It is clear that $W$ is an $\mathcal F$-maximal subgroup of $K\cap L\unlhd L$. Moreover, since $\mathfrak N^\pi$ is closed under taking subgroups, $L^{\mathfrak N^\pi}\le G^{\mathfrak N^\pi}\cap L\le K\cap L$.

    \item If $W\le X\le G$, $X\in \mathcal F$, then $X\cap K=W$.
    \smallskip

    We have that $W\le X\cap K\le K$, and $X\cap K\in \mathcal F$, because $X\cap K\unlhd X\in \mathcal F$. Since $W$ is $\mathcal F$-maximal in $K$ we deduce that $W=X\cap K$.

\item Whenever $U/W$ is $\mathfrak N^\pi$-Dsubnormal in $R^*/W\le L/W$, then $U=R^*_{\mathcal F}$ and $R^*\le N_L(U)$. In particular, this holds if $U/W\le R^*/W\in \mathfrak N^\pi$, $R^*/W\le L/W$. Also, if $H/W\le L/W$ and $R^*/W$ is $\mathfrak N^\pi$-Dnormal in $H/W$, then $H\le N_L(U)$. Moreover, $N_G(R^*)\le N_G(U)$.
   \smallskip

   By Lemma~\ref{lem2}(5) we deduce that $U$ is $\mathfrak N^\pi$-Dsubnormal in $R^*\le L$. Since  $U$ is $\mathcal F$-maximal in $L$, Proposition~\ref{cor1} implies that $U=R^*_{\mathcal F}$. If $U/W\le R^*/W\in \mathfrak N^\pi$, Lemma~\ref{lem2}(6) implies that $U/W$ is $\mathfrak N^\pi$-Dsubnormal in $R^*/W$.
   Obviously, if $R^*/W$ is $\mathfrak N^\pi$-Dnormal in $H/W$, then $U/W$ is $\mathfrak N^\pi$-Dsubnormal in $H/W$ and $H\le N_L(U)$, as above. Moreover, by Lemma~\ref{radical} it follows  that $N_G(R^*)\le N_G(R^*_{\mathcal F})=N_G(U)$.

        \item The following statements are pairwise equivalent:\begin{description}

\item (i) There exists $R\in \proj_{\mathfrak N^\pi}(L)$ such that $U= (RW)_{\mathcal F}$.

\item (ii) There exists $W\le R^*\le L$ such that $R^*/W\in \proj_{\mathfrak N^\pi}(L/W)$ and $U= R^*_{\mathcal F}$.

\item (iii)  There exists $W\le R^*\le L$ such that $U/W\le R^*/W\in \proj_{\mathfrak N^\pi}(L/W)$.
\smallskip

In this case, $\proj_{\mathfrak N^\pi}(N_L(U)/W)\subseteq \proj_{\mathfrak N^\pi}(L/W)$.
\end{description}
\smallskip

The equivalence (i) $\Leftrightarrow$ (ii) is clear by \cite[III. Proposition (3.7)]{DH}. The equivalence (ii) $\Leftrightarrow$ (iii) is a consequence of the fact 3. Finally, if we assume (ii), $R^*/W\in \proj_{\mathfrak N^\pi}(L/W)$ and $U= R^*_{\mathcal F}$, then $R^*\le N_L(U)$. By Lemma~\ref{teo1}, $R^*/W\in \proj_{\mathfrak N^\pi}(N_L(U)/W)$, and also
$\proj_{\mathfrak N^\pi}(N_L(U)/W)=\{(R^*/W)^x\mid x\in N_L(U)\}\subseteq \proj_{\mathfrak N^\pi}(L/W).$
        \end{enumerate}

By fact 2, $V\cap K=W$. We claim that $V/W\le Z_{\mathfrak N^\pi}(N_G(V)/W)$ the $\mathfrak N^\pi$-hypercentre of $N_G(V)/W$ (see \cite[IV. Notation and Definitions (6.8)]{DH}). Set $N=N_G(V)$. Since $N/(N\cap K)\cong NK/K\le G/K\in \mathfrak N^\pi$, we have that $N/(N\cap K)\in \mathfrak N^\pi$, since $\mathfrak N^\pi$ is closed under taking subgroups, and we can deduce  that $N$ acts $\mathfrak N^\pi$-hypercentrally on $N/(N\cap K)$. Then  $N$ acts $\mathfrak N^\pi$-hypercentrally on $V(N\cap K)/(N\cap K)$ which is $N$-isomorphic to $V/({V\cap N\cap K})=V/(K\cap V)=V/W$. It follows that $V/W\le Z_{\mathfrak N^\pi}(N/W)$, which proves the claim. Consequently, $V/W\le P^*/W\in \proj_{\mathfrak N^\pi}(N_G(V)/W)$ (see \cite[IV. Theorem (6.14)]{DH}).
\smallskip

We aim to prove that $V/W\le P^*/W\in \proj_{\mathfrak N^\pi}(G/W)$, which will conclude the proof by fact 4.

\smallskip

We prove next that $P^*/W$ is self-$\mathfrak N^\pi$-Dnormalizing in $G/W$; in particular we will have that $P^*/W$ is an $\mathfrak N^\pi$-maximal subgroup of $G/W$ by \cite[Corollary 3.10]{AADFP}. Assume that $P^*/W$ is $\mathfrak N^\pi$-Dnormal in $H/W\le G/W$. By fact 3, $H\le N_G(V)$. Since $P^*/W\in \proj_{\mathfrak N^\pi}(N_G(V)/W)$, Lemma~\ref{teo2} implies that $P^*/W=H/W$, and $P^*/W$ is self-$\mathfrak N^\pi$-Dnormalizing in $G/W$.
\smallskip

For any $H\le G$, we denote $\overline{H}=HW/W$. We distinguish next the following cases:
\begin{description}
\item Case 1: $\overline{G}=\overline{G}^{\mathfrak N^\pi}\,\overline{N_G(V)}$.
\item Case 2: $\overline{G}^{\mathfrak N^\pi}\,\overline{N_G(V)}<\overline{G}$.
\end{description}

\noindent
Case 1: $\overline{G}=\overline{G}^{\mathfrak N^\pi}\,\overline{N_G(V)}$.
\smallskip

In this case, $\G=\G^{\mathfrak N^\pi}\N^{\mathfrak N^\pi}\P^*=\G^{\mathfrak N^\pi}\P^*$, because $\mathfrak N^\pi$ is closed under taking subgroups and so $\N^{\mathfrak N^\pi}\le \G^{\mathfrak N^\pi}$. Then \begin{align*}\G/(\G^{\mathfrak N^\pi})'&=(\G^{\mathfrak N^\pi}/(\G^{\mathfrak N^\pi})')(\P^*(\G^{\mathfrak N^\pi})'/(\G^{\mathfrak N^\pi})')\\&= (\G/(\G^{\mathfrak N^\pi})')^{\mathfrak N^\pi}(\P^*(\G^{\mathfrak N^\pi})'/(\G^{\mathfrak N^\pi})').\end{align*} Let $\Q/(\G^{\mathfrak N^\pi})'$ be an $\mathfrak N^\pi$-maximal subgroup of $\G/(\G^{\mathfrak N^\pi})'$ such that $$\P^*(\G^{\mathfrak N^\pi})'/(\G^{\mathfrak N^\pi})'\le \Q/(\G^{\mathfrak N^\pi})'.$$ Then
$\Q/(\G^{\mathfrak N^\pi})'\in \proj_{\mathfrak N^\pi}(\G/(\G^{\mathfrak N^\pi})')$ by \cite[III. Lemma (3.14)]{DH}.

We consider the following two possibilities for $\Q$:
\smallskip

\underline{Case 1.1:} $\Q< \G$.
\smallskip

Let $W\le Q\le G$ such that $\Q=Q/W$. Then $W\le V\le P^*\le Q<G$. In particular, $P^*/W\le N_Q(V)/W\le N_G(V)/W$, which implies that $P^*/W\in \proj_{\mathfrak N^\pi}(N_Q(V)/W)$ because $P^*/W\in \proj_{\mathfrak N^\pi}(N_G(V)/W)=\cov_{\mathfrak N^\pi}(N_G(V)/W)$ by Lemma~\ref{teo1}. By inductive hypothesis (fact 1) and fact 4, we deduce that $\P^*=P^*/W\in \proj_{\mathfrak N^\pi}(Q/W)=\proj_{\mathfrak N^\pi}(\Q)$. But $\Q/(\G^{\mathfrak N^\pi})'\in \proj_{\mathfrak N^\pi}(\G/(\G^{\mathfrak N^\pi})')$, which implies that $\P^*\in \proj_{\mathfrak N^\pi}(\G)$, by \cite[III. Proposition (3.7)]{DH}, and we are done.
\smallskip

\underline{Case 1.2:} $\Q=\G$.
\smallskip

If $\Q=\G$, then $\G/(\G^{\mathfrak N^\pi})'\in \mathfrak N^\pi$ and so $\G^{\mathfrak N^\pi}=(\G^{\mathfrak N^\pi})'$.

Consequently,

$$\P^* \cap \G^{\mathfrak N^\pi}=\P^* \cap (\G^{\mathfrak N^\pi})'\le (\G^{\mathfrak N^\pi})'.$$

Assume that $\P^* \le \overline{X}< \G$. We may argue as above, in Case~1.1, to deduce by inductive hypothesis that $\P^*\in \proj_{\mathfrak N^\pi}(\overline{X})$.
Therefore, by Lemma~\ref{teo2}, it follows that  $$\P^* \cap \overline{X}^{\mathfrak N^\pi}\le (\overline{X}^{\mathfrak N^\pi})'.$$

Since $\P^*$ is self-$\mathfrak N^\pi$-Dnormalizing in $\G$, again Lemma~\ref{teo2} implies that $\P^*\in \proj_{\mathfrak N^\pi}(\G)$, which concludes the proof.
\smallskip

It remains to consider Case 2.
\smallskip

\noindent
Case 2: $\overline{G}^{\mathfrak N^\pi}\,\overline{N_G(V)}<\overline{G}$.
\smallskip

Since $\V \in \mathfrak N^\pi$, we can write $\V = \Vpi \times \Vpin$ where $\Vpi = O_{\pi}(\V)$ and $\Vpin = O_{\pi'}(\V)\in \mathfrak N$. Moreover, if $\G_{\pi}$ is a Hall $\pi$-subgroup of $\G$, then $\G^{\mathfrak N^\pi}\G_{\pi}\unlhd \G$, and we can form the subgroup $\G^{\mathfrak N^\pi}\G_{\pi}\Vpin$, which clearly contains $\V$. Assume that $\G^{\mathfrak N^\pi}\G_{\pi}\Vpin<\G$, and let $\Ls$ be a maximal subgroup of $\G$ containing $\G^{\mathfrak N^\pi}\G_{\pi}\Vpin$. Since $\G/\G^{\mathfrak N^\pi}\G_{\pi}\in \mathfrak N$, we have that $\Ls\unlhd \G$. By inductive hypothesis, $\V\le \R^*$ for some $\overline{R^*}\in \proj_{\mathfrak N^\pi}(\Ls)$. Hence, by the Frattini Argument (see \cite[A.(5.13)]{DH}) and fact 3,
$\G=\Ls N_{\G}(\R^*)=\Ls^{\mathfrak N^\pi}\R^*\N=\G^{\mathfrak N^\pi}\N$, which is not the considered case. Therefore, $\G^{\mathfrak N^\pi}\G_{\pi}\Vpin=\G$. In particular, $O^{\pi}(\G)=\G^{\mathfrak N^\pi}\Vpin\le \G^{\mathfrak N^\pi}\N<\G$, and $O^{\pi}(G)W=G^{\mathfrak N^\pi}V_{\pi'}W<G$, where $V_{\pi'}\le V$ such that $\Vpin=V_{\pi'}W/W$.

Denote $L=O^{\pi}(G)W=G^{\mathfrak N^\pi}V_{\pi'}W<G$. We prove next that $V_{\pi'}W$ is an $\mathcal F$-maximal subgroup of $L$. Since $V_{\pi'}W\unlhd V\in \mathcal F$, we have that $V_{\pi'}W\in \mathcal F$. Assume that $V_{\pi'}W\le  X\le L$ with $X\in \mathcal F$. By fact 2, $W=X\cap K$; in particular, $X\cap G^{\mathfrak N^\pi}\le W$. Hence $X=V_{\pi'}W(G^{\mathfrak N^\pi}\cap X)= V_{\pi'}W$. We deduce by inductive hypothesis (fact 1) and fact 3 that $\Vpin\le \R^*$ for some $\R^*\in \proj_{\mathfrak N^\pi}(\Ls)$, and $N_{\G}(\overline{R^*})\le N_{\G}(\Vpin)$. Consequently, $\G= \Ls N_{\G}(\R^*)=\G^{\mathfrak N^\pi} N_{\G}(\Vpin)$.

\smallskip

If $N_{\G}(\Vpin)<\G$, since $\V\le N_{\G}(\Vpin)$, by inductive hypothesis (fact 1) and fact 3 we have that $\V\le \R^*\in \proj_{\mathfrak N^\pi}(N_{\G}(\Vpin))$ and $\R^*\le \N$. Then
$\G=\G^{\mathfrak N^\pi} N_{\G}(\Vpin)=\G^{\mathfrak N^\pi}\R^*=\G^{\mathfrak N^\pi}\N $, which is not the case.
\smallskip

Assume finally that $\Vpin\unlhd \G$. Then $\Vpi\cong \V/\Vpin\le \G_{\pi}\Vpin/\Vpin$ for some  Hall $\pi$-subgroup $G_{\pi}$ of $G$. But then $\V/\Vpin$ is $\mathfrak N^\pi$-Dnormal in $\G_{\pi}\Vpin/\Vpin$, which implies by Lemma~\ref{lem2}(4) that $\V$ is $\mathfrak N^\pi$-Dnormal in $\G_{\pi}\Vpin$, and so $V=(G_{\pi}V_{\pi'}W)_{\mathcal F}$ by fact 3. Consequently, $\V\unlhd \G_{\pi}\Vpin$, and so $\G=\G^{\mathfrak N^\pi}\N$, the final contradiction.
\smallskip

(b) Let $G^*=\langle V,V_1\rangle$, $K^*=K\cap G^*\unlhd G^*$. Hence, $G^*/K^*\cong G^*K/K\le G/K\in \npi$, which implies that $G^*/K^*\in \npi$, since $\npi$ is closed under taking subgroups, and so $(G^*)^{\npi}\le K^*$. Also, $W$ is $\mathcal F$-maximal in $K^*$, because $W\le K\cap G^*=K^*$ and $W$ is $\mathcal F$-maximal in $K$. As in Part~(a),~fact~2, we can deduce that $K^*\cap V= K^*\cap V_1=W$. Consequently, $W\unlhd \langle V,V_1\rangle=G^*$. By Part~(a) there exist $\npi$-projectors $P$ and $Q$ of $G^*$ such that $V=(WP)_{\mathcal F}$ and $V_1=(WQ)_{\mathcal F}$. By Lemma~\ref{teo1} there exists $x\in G^*$ such that $Q^x=P$. Moreover, note that $G^*=(G^*)^{\npi}Q$, so that we can take $x\in (G^*)^{\npi}$. Consequently, by Lemma~\ref{radical}, it follows that
$$V_1^x= ((WQ)_{\mathcal F})^x= ((WQ)^x)_{\mathcal F}=(WQ^x)_{\mathcal F}=(WP)_{\mathcal F}=V,$$
with $x\in (G^*)^{\npi}$.\qed

\begin{teo}\label{teo*} If $\mathcal F$ is an $\npi$-Fitting set of a $\pi'$-soluble group $G$, then $G$ possesses exactly one conjugacy class of $\mathcal F$-injectors. Moreover, if $V$ and $V^*$ are $\mathcal F$-injectors of $G$, there exists $g\in G^{\npi}$ such that $(V^*)^g=V$.

Also, if $I$ is an $\mathcal F$-injector of $G$ and $N$ is an $\npi$-Dnormal subgroup of $G$, then $I\cap N$ is an $\mathcal F$-injector of $N$.
\end{teo}
\smallskip

{\it Proof.} We argue by induction on $|G|$. We may assume that $|G|\neq 1$ and that the result holds for all proper subgroups of $G$. Since $G$ is $\pi'$-soluble, $K=G^{\npi}$ is a normal proper subgroup of $G$. Let $W\in \inj_{\mathcal F}(K)$ and $V$ be an $\mathcal F$-maximal subgroup of $G$ containing $W$. We aim to prove first that $V\cap N\in \inj_{\mathcal F}(N)$ whenever $N$ is an $\npi$-Dnormal subgroup of $G$; in particular, $V\in \inj_{\mathcal F}(G)$ and $V$ satisfies the property stated in the last part of the statement. Let $M$ be a maximal $\npi$-Dnormal proper subgroup of $G$. It is enough to prove that $V\cap M\in \inj_{\mathcal F}(M)$. Note that $K=G^{\npi}\le M$, by Lemma~\ref{maxFDn}. Let $V_0\in \inj_{\mathcal F}(M)$. Then $V_0\cap K\in \inj_{\mathcal F}(K)$ and by inductive hypothesis we deduce that $W=(V_0\cap K)^g=V_0^g\cap K$ for some $g\in K$. We may replace $V_0$ by $V_0^g$, if necessary, and suppose that $V_0\cap K=W$. Let $V_1$ be an $\mathcal F$-maximal subgroup of $G$ such that $V_0\le V_1$. By Lemma~\ref{lema*}(b) and taking into account that $\npi$ is closed under taking subgroups, there exists $x\in \langle V,V_1\rangle^{\npi}\le G^{\npi}=K\le M$ such that $V_1^x=V$. Consequently, $V_0^x=V_0^x\cap M\le V_1^x\cap M=V\cap M$. By Lemma~\ref{lem3}(3), $V\cap M$ is $\npi$-Dnormal in $V\in \mathcal F$, so that $V\cap M\in \mathcal F$. But also $V_0^x\in \inj_{\mathcal F}(M^x)=\inj_{\mathcal F}(M)$, which implies that $V_0^x=V\cap M\in \inj_{\mathcal F}(M)$ as claimed.
\smallskip

We prove finally the conjugacy of $\mathcal F$-injectors, i.e. assume that $V^*\in \inj_{\mathcal F}(G)$ and prove that there exists $g\in G^{\npi}=K$ such that
$(V^*)^g=V$. It holds that $V^*\cap K\in \inj_{\mathcal F}(K)$. Then the inductive hypothesis implies that $(V^*\cap K)^k=W$ for some $k\in K^{\npi}\le K$. We consider now $V$ and $(V^*)^k$, which are $\mathcal F$-maximal subgroups of $G$ containing $W$. By Lemma~\ref{lema*}(b), there exists $t\in \langle V,(V^*)^k\rangle^{\npi}\le G^{\npi}$ such that $(V^*)^{kt}=V$, with $kt\in G^{\npi}$,  which concludes the proof.\qed
\medskip

As in \cite[IX.1]{DH} we state now corresponding concepts and results for general Fitting classes from the theory of Fitting sets.
\smallskip

In a natural way, we define a non-empty class $\F$ to be  an $\npi$-Fitting class if the following conditions are satisfied:
\begin{enumerate}
\item[(i)] If $G\in \F$ and $N$ is an $\npi$-Dnormal subgroup of $G$, then $N\in \F$.
\item[(ii)] If $M,N$ are $\npi$-Dnormal subgroups of $G=\langle M,N\rangle$ with $M,N\in \F$, then $G\in \F$.
\end{enumerate}

Now we have that an $\npi$-Fitting class is a  Fitting class, and also Fitting classes are exactly $\mathfrak N$-Fitting classes, for $\npi=\mathfrak N$ with $|\pi|\le 1$. As we show in the remarks below, if $|\pi|\ge 2$, not every Fitting class is an $\npi$-Fitting class, and $\npi$ is the smallest $\npi$-Fitting class of full characteristic (cf. \cite[IX. Theorem (1.9)]{DH}).

If $\F$ is an $\npi$-Fitting class and $G$ a group, then the trace of $\F$ in $G$, that is the set $\text{Tr}_{\F}(G)=\{H\le G\mid H\in \F\}$, is a $\npi$-Fitting set of $G$, and $\F$-injectors of $G$ are exactly $\text{Tr}_{\F}(G)$-injectors.

From Theorem~\ref{teo*} we can derive now the following result for $\npi$-Fitting classes and $\pi'$-soluble groups.

\begin{cor}\label{cor*} Let  $\F$ be an $\npi$-Fitting class and $G$ be a $\pi'$-soluble group, then $G$ possesses exactly one conjugacy class of $\F$-injectors. Moreover, if $V$ and $V^*$ are $\F$-injectors of $G$, there exists $g\in G^{\npi}$ such that $(V^*)^g=V$.

Also, if $I$ is an $\F$-injector of $G$ and $N$ is an $\npi$-Dnormal subgroup of $G$, then $I\cap N$ is an $\F$-injector of $N$.
\end{cor}
\begin{remarks} 1. The class $\npi$ is a particular case of the so-called lattice formations, which are classes of groups whose elements are direct product of Hall subgroups corresponding to pairwise disjoint sets of primes. With the same flavour as $\npi$-Fitting classes, though within the universe of finite soluble groups,  $\mathfrak L$-Fitting classes, for general lattice formations $\mathfrak L$ of soluble groups, were already defined in~\cite{AP04}.
\smallskip

2. If $\mathfrak F$ is an $\npi$-Fitting class of characteristic $\tau$, then $\npi\cap \mathfrak E_\tau \subseteq \mathfrak F$. In particular, if $|\pi|\ge 2$, then $\mathfrak N$ or also $\mathfrak N^m$, for any integer $m>1$, are never $\npi$-Fitting classes.
\smallskip

{\it Proof.} (\cite[Proposition 3.6]{AP04}) Assume that $G$ is a group of minimal order in $\npi\cap \mathfrak E_\tau \setminus \mathfrak F$. By Lemma~\ref{lem2}(6), a maximal subgroup of $G$ is $\npi$-Dnormal in $G$. By the choice of $G$, there is a unique maximal subgroup of $G$, which implies that $G$ is a cyclic $p$-group for some prime $p\in \tau$. Then $G\in \mathfrak F$ by  \cite[IX. Lemma (1.8)]{DH}), a contradiction.\qed
\smallskip

3. $\npi$ is an $\npi$-Fitting class.
\smallskip

{\it Proof.}  Since $\npi$ is closed under taking subgroups, condition (i) of the definition of $\npi$-Fitting class is satisfied.

Assume now that $M,N$ are $\npi$-Dnormal subgroups of $G=\langle M,N\rangle$ with $M,N\in \npi$, and we aim to prove that $G\in \npi$. For any $X\in \{M,N\}$, let $X=X_{\pi}\times X_{\pi'}$ with $X_{\pi}=O_{\pi}(X)$, $X_{\pi'}=O_{\pi'}(X)=O^{\pi}(X)\in \mathfrak N$. Since $X$ is $\npi$-Dnormal in $G$, by Proposition~\ref{pro1} and Lemma~\ref{lem3}(1), we have that $X_{\pi'}\unlhd G$, $M_{\pi'}N_{\pi'} \le  O^\pi(G)\le N_G(X)\le N_G(X_\pi)$ and $\langle X^G\rangle /X_{\pi'}$ is a $\pi$-group. Hence  $G=M_{\pi'}N_{\pi'}\langle M_\pi,N_\pi\rangle$ with  $M_{\pi'}N_{\pi'}$ a normal nilpotent $\pi'$-subgroup of $G$ and $[M_{\pi'}N_{\pi'},\langle M_\pi,N_\pi\rangle]=1$. Moreover $\langle M^G\rangle N_{\pi'}/M_{\pi'}N_{\pi'}$ and $\langle N^G\rangle M_{\pi'}/M_{\pi'}N_{\pi'}$ are $\pi$-groups, which implies that $G/M_{\pi'}N_{\pi'}$ is a $\pi$-group, and so $G=M_{\pi'}N_{\pi'}G_\pi$ with $G_\pi$ a Hall $\pi$-subgroup of $G$, by the Schur-Zassenhaus theorem. Moreover, there exist $x,y\in M_{\pi'}N_{\pi'}$ such that $M_\pi \le G_\pi^x$ and $N_\pi \le G_\pi^y$, which implies that $M_\pi = M_\pi^{x^{-1}}\le G_\pi$ and also $N_\pi = N_\pi^{y^{-1}}\le G_\pi$. Therefore, $\langle M_\pi,N_\pi\rangle\le G_\pi$, and then $\langle M_\pi,N_\pi\rangle= G_\pi$. It follows  that $G=M_{\pi'}N_{\pi'}\times G_\pi\in \npi$, which concludes the proof.\qed
\end{remarks}

\bigskip

\noindent
{\bf Acknowledgments.} Research supported by Proyectos PROMETEO/2017/ 057 from the Generalitat Valenciana (Valencian Community, Spain), and PGC2018-096872-B-I00 from the Ministerio de Ciencia, Innovaci\'on y Universidades, Spain, and FEDER, European Union. The fourth author acknowledges with thanks the financial support of the Universitat de Val\`encia as research visitor (Programa Propi d'Ajudes a la Investigaci\'o de la Universitat de Val\`encia, Subprograma Atracci\'o de Talent de VLC-Campus, Estades d'investigadors convidats (2019)).

\bigskip

\noindent Milagros Arroyo-Jord\'a, Paz Arroyo-Jord\'a\\
Escuela T\'{e}cnica Superior de Ingenieros Industriales\\
Instituto Universitario de Matem\'{a}tica Pura y  Aplicada IUMPA\\
Universitat Polit\`ecnica de Val\`encia,\\
 Camino de Vera, s/n,  46022 Valencia, Spain\\
E-mail: marroyo@mat.upv.es, parroyo@mat.upv.es\\
 \\
Rex Dark\\
School of Mathematics, Statistics and Applied Mathematics,\\
National University of Ireland, University Road, Galway,
  Ireland\\
E-mail: rex.dark@nuigalway.ie\\
 \\
 Arnold D.~Feldman\\
Franklin and Marshall College, Lancaster, PA 17604-3003,
  U.S.A.\\
E-mail: afeldman@fandm.edu\\
 \\
Mar\'{\i}a Dolores P\'erez-Ramos\\
Departament de Matem\`{a}tiques, Universitat de Val\`{e}ncia,\\
C/ Doctor Moliner 50, 46100 Burjassot
(Val\`{e}ncia), Spain\\
E-mail: Dolores.Perez@uv.es

\end{document}